\documentclass[a4paper,11pt]{article}
\usepackage{amsmath,amssymb,epic,eepic,eepicemu,graphicx}
\usepackage[usenames]{color}
\usepackage[english]{babel}

\title{A conformal energy for simplicial surfaces}
\author{ Alexander I.\,Bobenko
\thanks{Partly supported by the DFG Research
Center ``Mathematics for key technologies'' (FZT 86) in Berlin}}
\date{\today}

\newtheorem{theorem}{Theorem}
\newtheorem{lemma}[theorem]{Lemma}
\newtheorem{proposition}[theorem]{Proposition}
\newtheorem{corollary}[theorem]{Corollary}
\newtheorem{definition}[theorem]{Definition}
\newtheorem{conjecture}[theorem]{Conjecture}

\begin{document}
\maketitle
\begin{center}
 Institut f\"ur Mathematik,
 %Fak. 2,
 Technische Universit\"at Berlin,\\
 Strasse~des 17.~Juni 136, 10623 Berlin, Germany\\
 \medskip e-mail:{\tt bobenko@math.tu-berlin.de}
\end{center}

%--------------------------------------------------------------------------
\begin{abstract}
A new functional for simplicial surfaces is suggested. It is
invariant with respect to M\"obius transformations and is a
discrete analogue of the Willmore functional. Minima of this
functional are investigated. As an application a bending energy
for discrete thin-shells is derived.
\end{abstract}

 \eject

%--------------------------------------------------------------------------
\section{Introduction}
\label{s.intro}

In the variational description of surfaces the following
functionals are of primary importance:
\begin{itemize}
    \item The area ${\cal A}=\int dA$, where $dA$ is the area
    element, is preserved by isometries.
    \item The total Gaussian curvature ${\cal G}=\int K dA$,
    where $K$ is the Gaussian curvature, is topological invariant.
    \item The total mean curvature ${\cal M}=\int H dA$, where $H$ is the
    mean curvature, depends on external geometry of the surface.
    \item The Willmore energy ${\cal W}=\int H^2 dA$ is invariant with
    respect to M\"obius transformations.
\end{itemize}

Geometric discretizations of the first three functionals for
simplicial surfaces are well known. For the area functional it is
obvious. The local Gaussian curvature at a vertex $v$ is defined
as the angle defect
$$ G(v)=2\pi - \sum_i \alpha_i, $$
 where
$\alpha_i$ are the angles of all triangles (see Fig.
\ref{f.alpha-beta}) at vertex $v$. The total Gaussian curvature is
a sum over all vertices $G=\sum_v G(v)$. The local mean curvature
at an edge $e$ is defined as
$$ M(e)=l\theta, $$
where $l$ is the length of the edge and $\theta$ is the angle
between the normals to the adjacent faces at $e$ (see Fig.
\ref{f.bending}). The total mean curvature is the sum over all
edges $M=\sum_e M(e)$. These discrete functionals possess the
geometric symmetries of the smooth functionals mentioned above.

Until recently a geometric discretization of the Willmore
functional was missing. In this paper we introduce a M\"obius
invariant energy for simplicial surfaces and show that it should
be treated as a discrete Willmore energy.
%--------------------------------------------------------------------------
\section{Conformal Energy}
\label{s.energy}

Let $S$ be a simplicial surface in 3-dimensional Euclidean space
with the set of vertices $V$, edges $E$ and (triangular) faces
$F$. We define a conformal energy for simplicial surfaces using
circumcircles of their faces. Each (internal) edge $e\in E$ is
incident to two triangles. Consistent orientation of the triangles
naturally induces an orientation of the corresponding
circumcircles. Let $\beta(e)$ be the external intersection angle
of the circumcircles of the triangles sharing $e$, which is the
angle between the tangent vectors of the oriented circumcircles.
Now we are ready to define a new energy functional for simplicial
surfaces.

\begin{definition}
Local conformal (discrete Willmore) energy at a vertex $v$ is
given by the sum over all incident edges
$$
W(v)=\sum_{e\ni v} \beta(e)-2\pi.
$$
The conformal (discrete Willmore) energy of the simplicial surface
is the sum over all vertices
$$
W(S)=\dfrac{1}{2}\sum_{v\in V}W(v)=\sum_{e\in E} \beta(e)-\pi\mid
V \mid,
$$
where $\mid V\mid$ is the number of vertices of $S$.
\end{definition}

\begin{figure}[h]
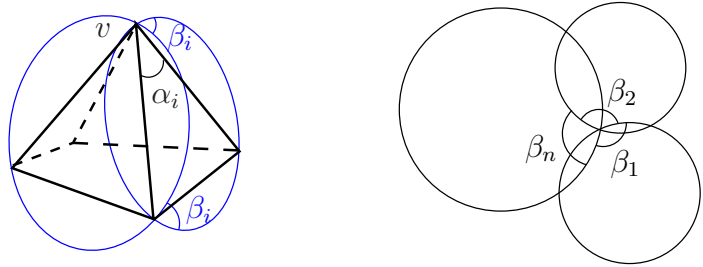

\begin{center}
\parbox[c]{0.4\textwidth}{\input{dwill2.pstex_t}}
%\hfill
\parbox[c]{0.3\textwidth}{\input{circles.pstex_t}}
\end{center}
\caption{Definition of conformal (discrete Willmore) energy
\label{f.dEnergy}}
\end{figure}

Fig.\ref{f.dEnergy} presents two neighboring circles with their
external intersection angle $\beta_i$ as well as a view ``from the
top'' at a vertex $v$ showing all $n$ circumcircles passing
through $v$ with the corresponding intersection angles
$\beta_1,\ldots,\beta_n$. For simplicity we will consider only
simplicial surfaces without boundary.

Note that the energy we defined is obviously invariant with
respect to M\"obius transformations. This invariance is an
important property of the classical Willmore energy defined for
smooth surfaces (see below).

Note also that the changing of the orientation of both circles
preserves the angle $\beta(e)$. This shows that the energy $W(S)$
is well defined for non-oriented simplicial surfaces as well.

The star $S(v)$ of the vertex $v$ is the subcomplex of $S$
comprised by the triangles incident with $v$. The vertices of
$S(v)$ are $v$ and all its neighbors. We call $S(v)$ convex if for
any its face $f\in F(S(v))$ the star $S(v)$ lies to one side of
the plane of $F$ and strictly convex if the intersection of $S(v)$
with the plane of $f$ is $f$ itself.

\begin{proposition} \label{p.non-negative}
The conformal energy is non-negative
$$W(v)\ge 0,$$
and vanishes if and only if the star $S(v)$ is convex and all its
vertices lie on a common sphere.
\end{proposition}

\begin{figure}[h]
\begin{center}
\parbox[c]{0.4\textwidth}{\input{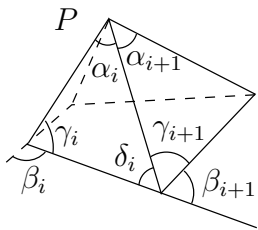}}
\end{center}
\caption{Proof of Lemma \ref{l.alpha-beta} \label{f.alpha-beta}}
\end{figure}

The proof of this proposition is based on the following elementary
lemma.
\begin{lemma} \label{l.alpha-beta}
Let $\cal P$ be a (not necessarily planar) $n$-gon with external
angles $\beta_i$. Choose a point $P$ and connect it to all
vertices of $\cal P$. Let $\alpha_i$ be the angles (see
Fig.\ref{f.alpha-beta}) of the triangles at the tip $P$ of the
obtained pyramid. Then
$$
\sum_{i=1}^n \beta_i\ge \sum_{i=1}^n \alpha_i,
$$
and the equality holds if and only if $\cal P$ is planar and
convex\footnote{The obtained pyramid is convex in this case. Note
that we distinguish convex and strictly convex polygons (and
pyramides). Some of the angles $\beta_i$ of a convex polygon may
vanish. The corresponding cite-triangles of the pyramid lie in one
plane.}
 and the vertex
$P$ lies inside $\cal P$.
\end{lemma}

{\it Proof}. Let us denote by $\gamma_i$ and $\delta_i$ the angles
of the cite-triangles at the vertices of $\cal P$ (see
Fig.\ref{f.alpha-beta}). The claim of Lemma \ref{l.alpha-beta}
follows from adding over all $i=1,\ldots,n$ two obvious
(in)equalities
\begin{eqnarray*}
\beta_{i+1}&\ge& \pi -(\gamma_{i+1}+\delta_i)\\
\pi -(\gamma_{i}+\delta_i)&=&\alpha_i.
\end{eqnarray*}
All inequalities become equalities only in the case when $\cal P$
is planar, convex and contains $P$. Note that some of the external
angles $\beta_i$ may vanish. As a corollary we obtain a polygonal
version of Fenchel's theorem \cite{Fenchel}.

\begin{corollary}  \label{c.non-negative}
$$
\sum_{i=1}^n \beta_i\ge 2\pi.
$$
\end{corollary}
{\it Proof.} For a given $\cal P$ choose the point $P$ varying on
a straight line encircled by $\cal P$. There always exist points
$P$ so that the star at $P$ is not strictly convex, and thus $\sum
\alpha_i \ge 2\pi$.

{\it Proof of Proposition \ref{p.non-negative}}. The claim of
Proposition \ref{p.non-negative} is invariant with respect to
M\"obius transformations. Applying a M\"obius transformation $M$
which maps the vertex $v$ to infinity, $M(v)=\infty$, one observes
that all the circles passing through $v$ become straight lines and
we arrive at the geometry shown in Fig.\ref{f.alpha-beta} with
$P=M(\infty)$. Now the claim follows immediately from Corollary
\ref{c.non-negative}.

\begin{theorem}
Let $S$ be a simplicial surface. Then
$$W(S)\ge 0,$$
and the equality holds if and only if $S$ is a (part of a) convex
polyhedron inscribed in a sphere.
\end{theorem}
{\it Proof.} Only the second statement needs to be proven. Due to
Proposition \ref{p.non-negative}  the equality $W(S)=0$ implies
that all vertices and edges of $S$ are convex (but not necessarily
strictly convex). Deleting the edges which separate triangles
lying in one plane one obtains a polyhedral surface $S_P$ with
circular faces and all strictly convex vertices and edges.
Proposition \ref{p.non-negative} implies that for every vertex $v$
there exists a sphere $S_v$ with all vertices of the star $S(v)$
lying on it. For any edge $(v_1,v_2)$ of $S_P$ two neighboring
spheres $S_{v_1}$ and $S_{v_2}$ share two different circles of
their common faces. This implies $S_{v_1}=S_{v_2}$ and finally the
coincidence of all the spheres $S_v$.

Discrete conformal energy $W$ defined above is a discrete analogue
of the Willmore energy for smooth surfaces \cite{Willmore}
$$
{\cal W(S)}=\frac{1}{4}\int_{\cal S} (k_1-k_2)^2 dA=\int_{\cal S}
H^2 dA- \int_{\cal S} K dA.
$$
Here $dA$ is the area element, $k_1, k_2$ the principal
curvatures, $H=\frac{1}{2}(k_1+k_2)$ the mean curvature, $K=k_1
k_2$ the Gaussian curvature of the surface. Here we prefer a
definition for $\cal W$ with a M\"obius invariant integrand. It
differs from the one in the introduction by a topological
invariant.

%$\int K=2\pi \chi({\cal S})$ topological invariant
Let us mention two important properties of the Willmore energy:
\begin{itemize}
    \item  ${\cal W(S)}\ge 0$ and ${\cal W(S)}=0$ if and only if
    $\cal S$ is the round sphere.
    \item  $\cal W(S)$ (and the integrand $(k_1-k_2)^2 dA$) is M\"obius
    invariant \cite{Blaschke,Willmore}
\end{itemize}
Whereas the first claim almost immediately follows from the
definition, the second one is a non-trivial property. We have
shown that the same properties hold for the discrete energy $W$;
in the discrete case the M\"obius invariance is built into the
definition and the non-negativity of the energy is non-trivial.

In the same way one can define conformal
    (Willmore) energy for simplicial surfaces in
    Euclidean spaces of higher dimension and space forms.

The discrete conformal energy is well defined for
    polyhedral surfaces with circular faces (not necessarily
    simplicial).

\section{Computation of the Energy}

\begin{figure}[h]
\begin{center}
\parbox[c]{0.4\textwidth}{\input{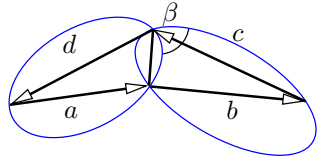}}
\end{center}
\caption{Formula for the angle between circumcircles
\label{f.formula}}
\end{figure}

Consider two triangles with a common edge. Let $a, b, c, d\in
{{\mathbb R}^3}$ be their other edges oriented as in
Fig.\ref{f.formula}. Identifying vectors in ${\mathbb R}^3$ with
imaginary quaternions ${\rm Im}\ \mathbb{H}$ one obtaines for the
quaternionic product
\begin{equation}                \label{product}
ab=-<a,b>+a\times b,
\end{equation}
where $<a,b>$ and $a\times b$ are the scalar and vector products
in ${\mathbb R}^3$.

\begin{proposition}
The external angle $\beta\in [0,\pi ]$ between the circumcircles
of the triangles in Fig.\ref{f.formula} is given by one of the
equivalent formulas:
\begin{eqnarray*}
\cos(\beta)&=&-\dfrac{{\rm Re}\ q}{\mid q\mid}=
            -\dfrac{{\rm Re}\ abcd}{\mid abcd\mid}= \\
           &=&\dfrac{<a,c><b,d>-<a,b><c,d>-<b,c><d,a>}
           {\mid a\mid \mid b\mid \mid c\mid \mid d\mid}.
\end{eqnarray*}
Here $q=ab^{-1}cd^{-1}$ is the cross-ratio of the quadrilateral.
\end{proposition}
{\it Proof.} Since ${\rm Re}\ q$, $\mid q\mid$ and $\beta$ are
M\"obius invariant it is enough to prove the first formula for the
planar case $a,b,c,d\in {\mathbb C}$, mapping all four vertices to
a plane by a M\"obius transformation. In this case $q$ becomes the
classical complex cross-ratio. Considering the arguments
$a,b,c,d\in {\mathbb C}$ one easily arrives at $\beta=\pi-\arg q$.
The second representation follows from the identity
$b^{-1}=-b/\mid b\mid$ for imaginary quaternions. Finally applying
(\ref{product}) we obtain
\begin{eqnarray*}
{\rm Re}\ abcd=<a,b><c,d>-<a\times b,c\times d>=\\
<a,b><c,d>+<b,c><d,a>-<a,c><b,d>. \end{eqnarray*}

\section{Minimizing Discrete Conformal Energy}

Similarly to the smooth Willmore functional $\cal W$, minimizing
the discrete conformal energy $W$ makes the surface as round as
possible.

Let us denote by $\bf S$ the combinatorial data of $S$. The
simplicial surface $S$ is called a geometric realization of the
abstract simplicial surface $\bf S$.

\begin{definition}
Critical points of $W(S)$ are called simplicial Willmore surfaces.
The conformal (Willmore) energy of an abstract simplicial surface
is the infimum over all geometric realizations
$$
W({\bf S})=\inf_{S\in {\bf S}} W(S).
$$
\end{definition}

\begin{figure}[h]
 \begin{center}
\begin{tabular}{c c c c c}
\includegraphics[width=0.25\textwidth]{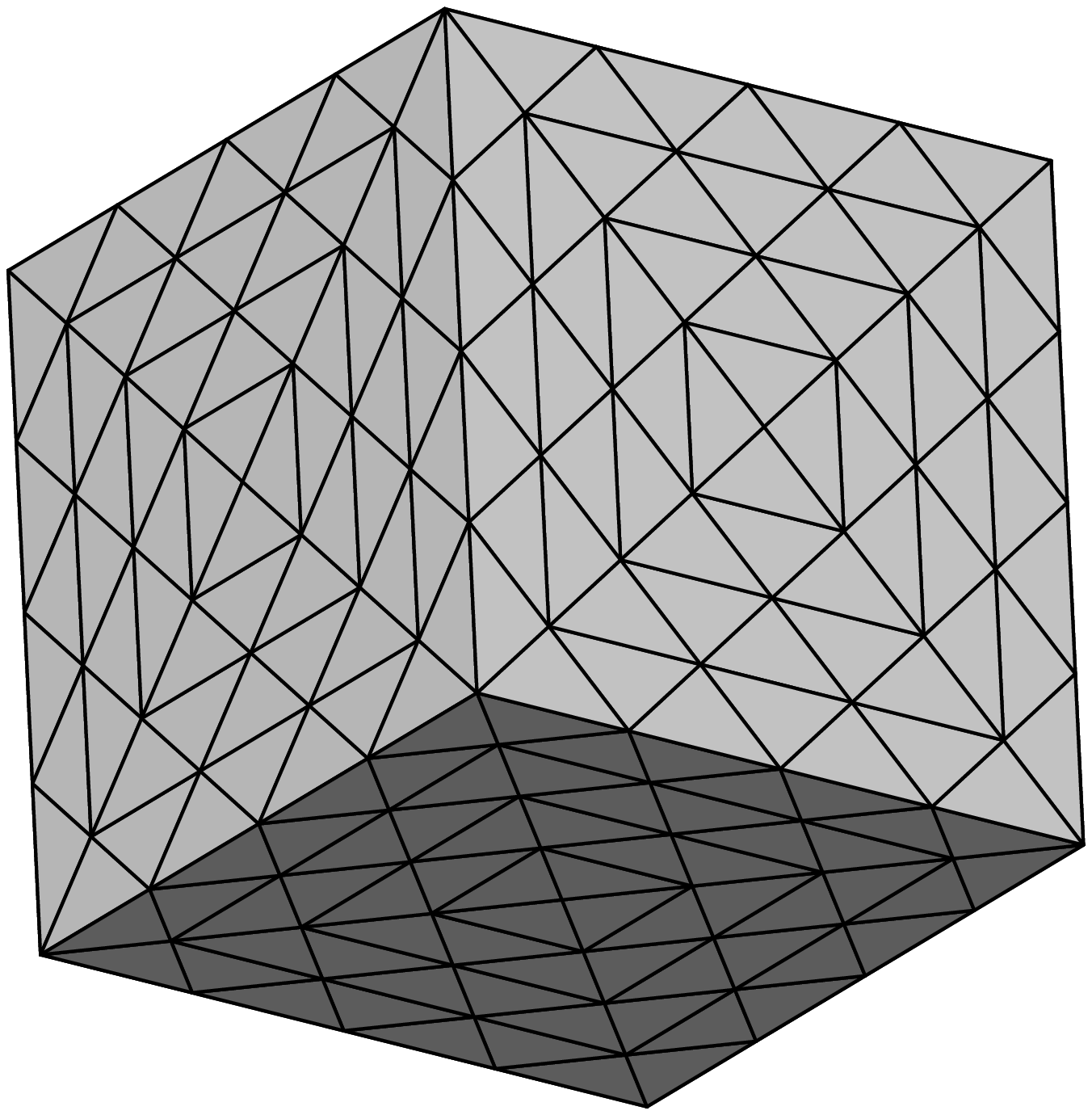}&\hfill&
\includegraphics[width=0.25\textwidth]{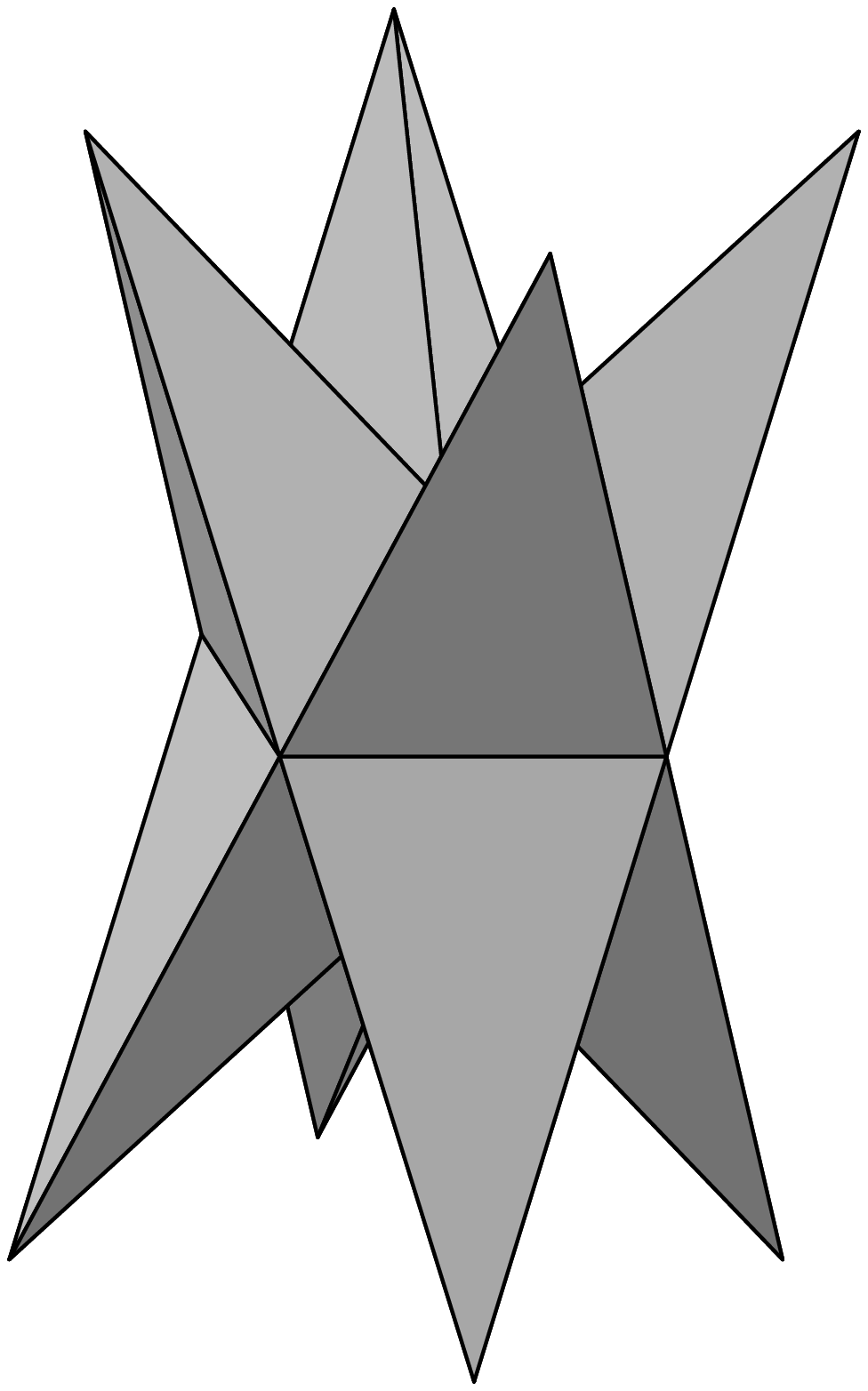}&\hfill&
\includegraphics[width=0.25\textwidth]{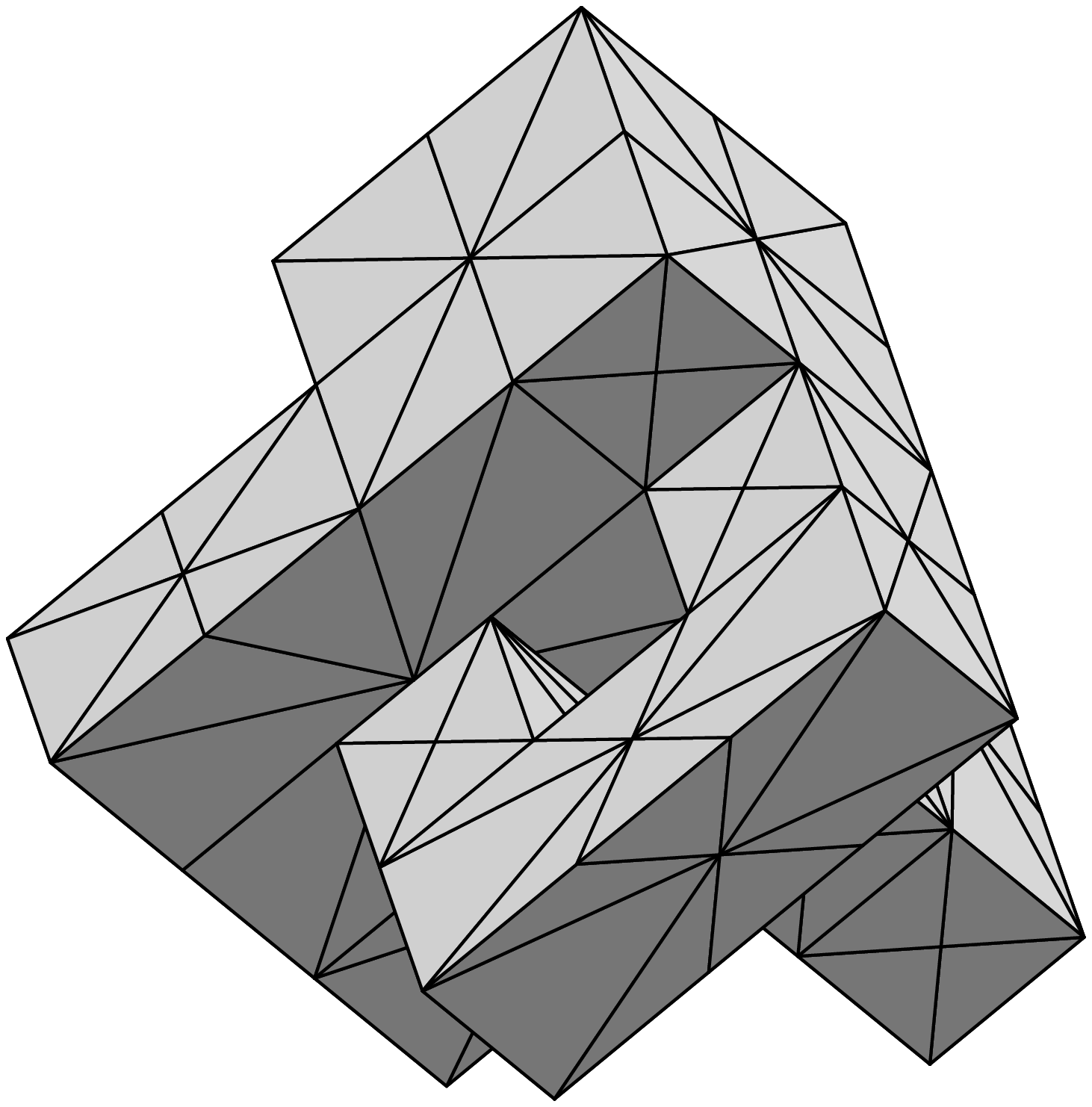}\\
\includegraphics[width=0.25\textwidth]{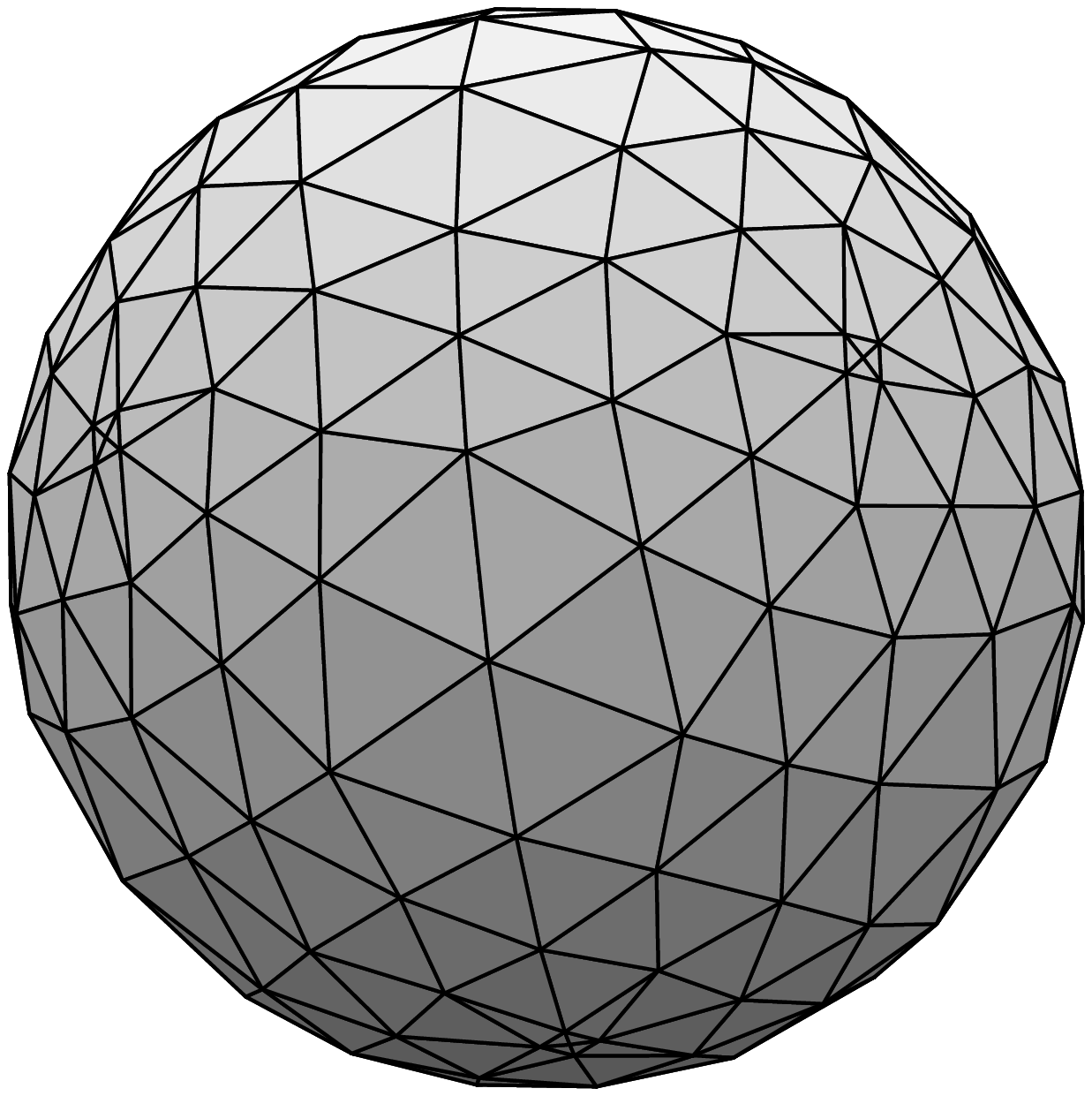}&\hfill&
\includegraphics[width=0.25\textwidth]{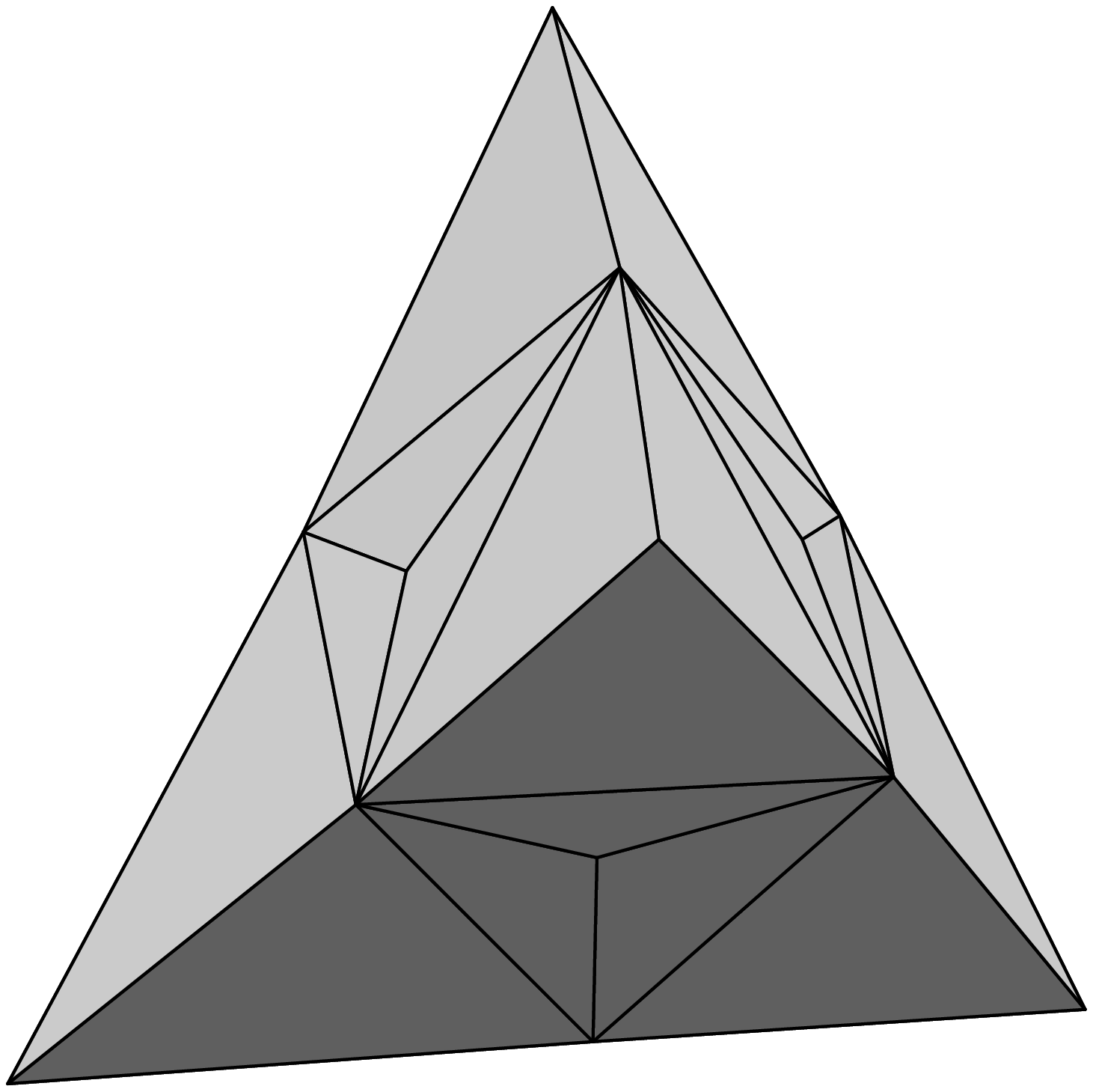}&\hfill&
\includegraphics[width=0.25\textwidth]{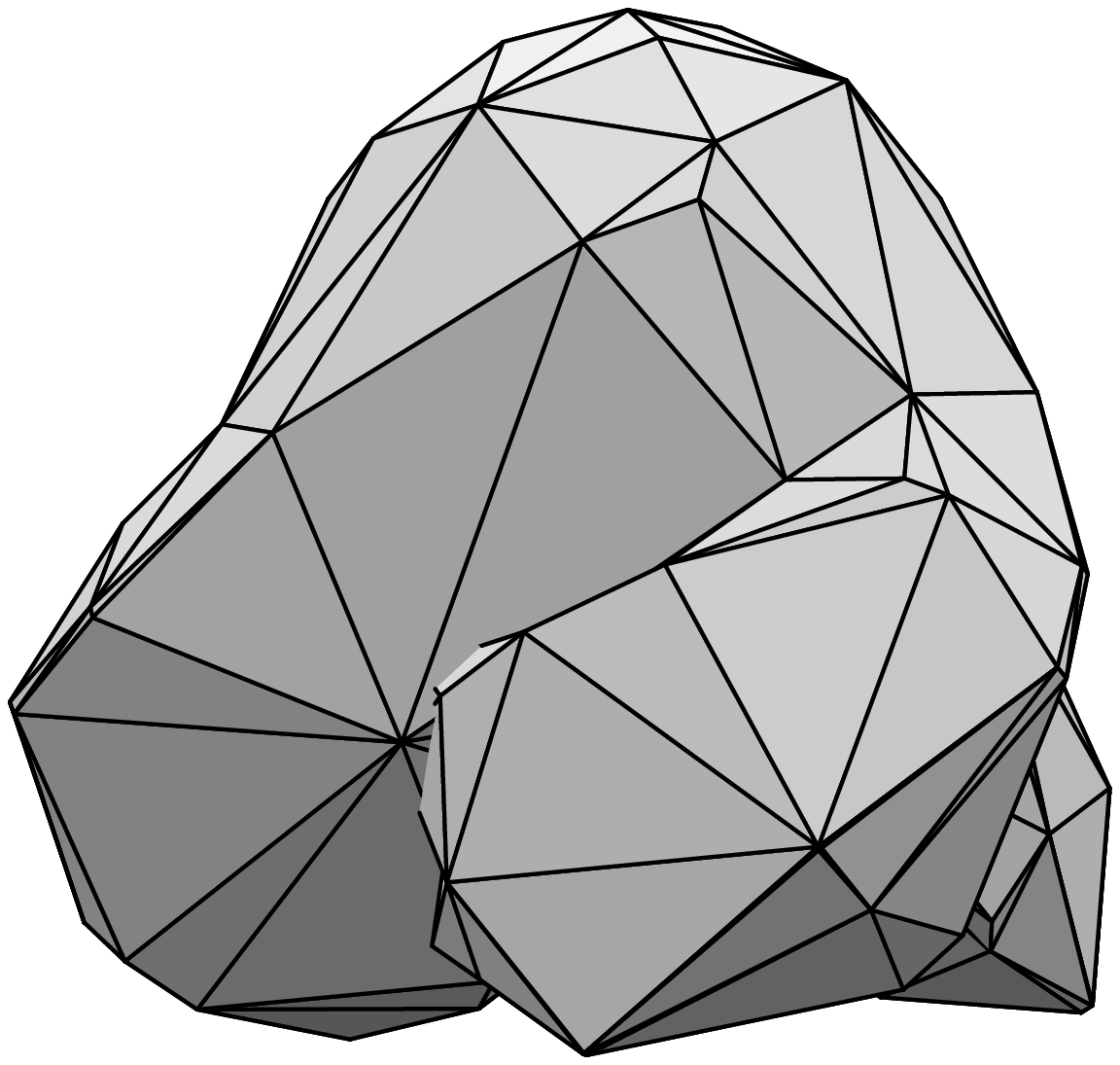}
%\\
%Triangulation:  &\hfill& Triangulation: &\hfill& discrete Boy
%surface\\
%``spherical'', $W=0$ & \hfill & ``non-spherical'', $W>0$ & \hfill
%& (projective plane)
\end{tabular}
\end{center}
\caption{Discrete Willmore spheres of ``inscribable'' ($W=0$) and
``non-inscribable'' ($W>0$) type and discrete Boy surface
\label{f.minimizing}}
\end{figure}

Kevin Bauer implemented the proposed conformal functional with the
Brakke's evolver \cite{Brakke} and did some minimization numerical
experiments. Examples of those are presented in
Fig.\ref{f.minimizing}. In the first and the second lines the
initial configurations and the corresponding Willmore surfaces
minimizing the conformal energy respectively are shown. Let us
call the gradient flow of the energy $W$ the discrete Willmore
flow. By this flow the energy of the first simplicial sphere
decreases to zero and the surface evolves into a convex polyhedron
with all the vertices lying on a sphere. The abstract simplicial
surface of the central example is different and we obtain a
simplicial Willmore sphere with positive conformal energy. The
third example is a simplicial projective plane. The initial
configuration is made from squares divided into triangles (cf.
\cite{Petit}). We see that the minimum is close to the smooth Boy
surface known to minimize \cite{KarcherP} the Willmore energy for
projective planes.

The minimization of the conformal energy for simplicial spheres is
related to a classical result of Steinitz \cite{Steinitz} who has
shown that there exist abstract simplicial 3-polytopes without
geometric realizations all vertices of which belong to a sphere.
We call these combinatorial types non-inscribable.

The non-inscribable examples of Steinitz are constructed as
follows \cite{Gruenbaum}. Let $\bf S$ be an abstract simplicial
sphere with vertices colored in black and white. Denote the sets
of white and black vertices by $V_w$ and $V_b$ respectively,
$V=V_w\cup V_b$. Assume that the number of black vertices does not
exceed the number of white vertices, $\mid V_w\mid \ge \mid
V_b\mid $, and there are no edges connecting two white vertices
and there are edges connecting black vertices. It is easy to see
that $\bf S$ with these properties cannot be inscribed in a
sphere. Indeed, assume that we have constructed such inscribed
convex polyhedron. Then the equality of the intersection angles at
both ends of an edge (see left Fig.\ref{f.dEnergy}) implies
$$
2\pi \mid V_b\mid \ge  \sum_{e\in E} \beta(e)\ge 2\pi \mid
V_w\mid,
$$
and the equalities hold only if all edges connect vertices of
different color. The obtained contradiction to the assumed
inequality implies the claim.

\begin{figure}[h]
 \begin{center}
\includegraphics[width=0.35\textwidth]{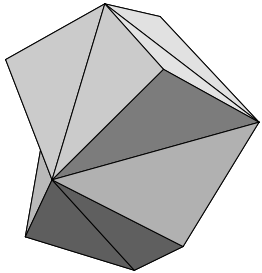}
\end{center}
\caption{A discrete Willmore sphere of ``non-inscribable'' type
with 11 vertices and $W=2\pi$ \label{f.11}}
\end{figure}

To construct abstract polyhedra with $\mid V_w\mid \ge \mid
V_b\mid $ and with some edges connecting black points, take a
polyhedron $\bf P$ whose number of vertices does not exceed the
number of faces $\mid \hat{F}\mid \ge \mid \hat{V}\mid $. Color
all the vertices in black, add white vertices at the faces and
connect them to all black vertices of a face. We obtain a
polyhedron with black (original) edges and $\mid V_w\mid=\mid
\hat{F}\mid\ge\mid V_b\mid=\mid \hat{V}\mid$. The example with
minimal possible number of vertices $\mid V\mid =11$ is shown in
Fig.\ref{f.11}. The starting polyhedron $\bf P$ here are two
tetrahedra identified along a common face: $\hat{F}=6, \hat{V}=5$.

Hodgson, Rivin and Smith \cite{HodgsonRS} found  a
characterization of inscribable combinatorial types, based on a
transfer to the Klein model of hyperbolic 3-space. It is not clear
whether there exist non-inscribable examples of non-Steinitz type.

Numerical experiments lead us to the following
\begin{conjecture}
The conformal energy of simplicial Willmore spheres is quantized
$$
W=2\pi N,\quad N\in {\mathbb N}.
$$
\end{conjecture}
Note that this claim belongs to differential geometry of discrete
surfaces. It would be interesting to find a (combinatorial)
meaning of the integer $N$.  Compare also with the famous
classification of smooth Willmore spheres by Bryant \cite{Bryant},
who has shown that the energy of Willmore spheres is quantized
${\cal W}= 4\pi N, N\in{\mathbb N}$.

The discrete Willmore energy is defined for the ambient spaces
($\mathbb{R}^n$ or $S^n$) of any dimension. This leads to
combinatorial Willmore energies
$$
W_n(S)=\inf_{S\in {\bf S}} W(S), \qquad S\subset S^n,
$$
where the infimum is taken over all realizations in the
$n-$dimensional sphere. Obviously these numbers build a
non-increasing sequence $W_n(S)\ge W_{n+1}(S)$ which becomes
constant for sufficiently large $n$.

Complete understanding of non-inscribable simplicial spheres is an
interesting mathematical problem. However the phenomenon of
existence of such spheres might be seen as a problem in using of
the conformal functional for applications in computer graphics,
such as fairing of surfaces. Fortunately the problem disappears
just after one refinement step: all simplicial spheres become
inscribable. Let $\bf S$ be an abstract simplicial sphere. Define
its refinement $\bf S_R$ as follows: split every edge of $\bf S$
in two by putting additional vertices and connect these new
vertices sharing a face of $\bf S$ by additional edges.

\begin{proposition}\label{p.refine}
The refined simplicial sphere $\bf S_R$ is inscribable, and thus
$W({\bf S_R})=0$.
\end{proposition}

{\it Proof.} Koebe's theorem (see, for example, \cite{Ziegler},
\cite{BobenkoS}) claims that every abstract simplicial sphere $\bf
S$ can be realized as a convex polyhedron $S$ all edges of which
touch a common sphere $S^2$. Starting with this realization $S$ it
is easy to construct a geometric realization $S_R$ of the
refinement $\bf S_R$ inscribed in $S^2$. Indeed choose the
touching points of the edges of $S$ with $S^2$ as additional
vertices of $S_R$ and project the original vertices of $S$ (which
lie outside of the sphere $S^2$) to $S^2$. One obtaines a convex
simplicial polyhedron $S_R$ inscribed in $S^2$.

Another interesting variational problem with the conformal energy
is the optimization of triangulations of a given simplicial
surface. Here one fixes the vertices and chooses an equivalent
triangulation (abstract simplicial surface $\bf S$) minimizing the
conformal functional. The minimum
$$
W(V)=\min_{{\bf S}\ni S} W(S)
$$
yields ``an optimal'' triangulation for a given vertex data. In
the case of $S^2$ this optimal triangulation is classical.

\begin{proposition}
Let $S$ be a simplicial surface with all vertices $V$ on a two
dimensional sphere $S^2$. Then $W(S)=0$ if and only if it is the
Delaunay triangulation on the sphere, i.e. $S$ is the boundary of
the convex hull of $V$.
\end{proposition}

In differential geometric applications like numerical minimizing
the Willmore energy of smooth surfaces (cf. \cite{HsuKS}) it is
not natural to preserve the triangulation by minimizing the
energy, and one should also change the combinatorial type
decreasing the energy.

Discrete conformal energy $W$ is not only a discrete analogue of
the Willmore energy. One can show that it approximates the smooth
Willmore energy although the smooth limit is very sensitive to the
refinement method and should be chosen in a very special way. A
computation which will be published elsewhere shows that if one
chooses the vertices of a curvature line net of a smooth surface
$\cal S$ for the vertices of $S$ and triangularizes it then $W(S)$
converges to $\cal W(S)$ by natural refinement. On the other hand
the infinitesimal equilateral triangular lattice gives in the
limit 3/2 times larger energy. Possibly the minimization of the
discrete Willmore energy with the vertices on the smooth surface
could be used for computation of the curvature line net. We are
going to investigate this interesting and complicated phenomenon.

\section{Bending of Simplicial Surfaces}

An accurate model for bending of discrete surfaces is important
for modelling in virtual reality.

Let ${\cal S}_0$ be a thin shell and ${\cal S}$ its deformation.
The bending energy of smooth thin shells is given by the integral
\cite{GrinspunHDS}
$$
E=\int (H-H_0)^2 dA,
$$
where $H_0$ and $H$ are the mean curvatures of the original and
deformed surface respectively. For $H_0=0$ it reduces to the
Willmore energy.

\begin{figure}[h]
\begin{center}
\parbox[c]{0.4\textwidth}{\input{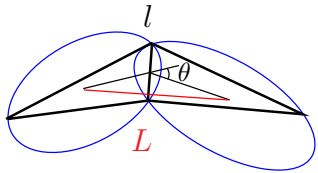}}
\end{center}
\caption{To definition of the bending energy for simplicial
surfaces \label{f.bending}}
\end{figure}

To derive the bending energy for simplicial surfaces let us
consider the limit of fine triangulation, i.e. of small angles
between the normals of neighboring triangles. Consider an
isometric deformation of two adjacent triangles. Let $\theta$ be
the complement of the dihedral angle of the edge $e$, or,
equivalently, the angle between the normals of these triangles
(see Fig.\ref{f.bending}) and $\beta(\theta)$ the external
intersection angle between the circumcircles of the triangles (see
Fig.\ref{f.dEnergy}) as a function of $\theta$.

\begin{proposition}  \label{p.bending}
Assume that the circumcenters of the circumcircles of two adjacent
triangles do not coincide. Then in the limit of small angles
$\theta\to 0$ the angle $\beta$ between the circles behaves as
follows:
$$
\beta(\theta)=\beta(0)+\dfrac{l}{L}\theta^2+o(\theta^3).
$$
Here $l$ is the length of the edge and $L\ne 0$ is the distance
between the centers of the circles.
\end{proposition}

This proposition and our definition of conformal energy for
simplicial surfaces motivate to suggest
$$
E=\sum_{e\in E}\dfrac{l}{L}\theta^2
$$
for the bending energy of discrete thin-shells.

In \cite{BridsonMF, GrinspunHDS} similar representations for the
bending energy of simplicial surfaces were found empirically. They
were demonstrated to give convincing simulations and good
comparison with real processes. In \cite{GrinspunHDS} the distance
between the barycenters is used for $L$ in the energy expression
but possible numerical advantages in using circumcenters are
indicated.

Using the Willmore energy and Willmore flow is a hot topic in
computer graphics. Applications include fairing of surfaces and
surface restoration. We hope that our conformal energy will be
useful for these applications and plan to work on them.

%\section{Discrete Minimal Surfaces}

{\bf Acknowledgements}. I would like to thank Ulrich Pinkall for
the discussion in which the idea of the discrete Willmore
functional was born. I am also grateful to G\"unter Ziegler, Peter
Schr\"oder, Boris Springborn, Yuri Suris and Ekkerhard Tjaden for
useful discussions and to Kevin Bauer for making numerical
experiments with the conformal energy.

%--------------------------------------------------------------------------

\end{document}